\documentclass{amsart}
\usepackage[T1]{fontenc}
\usepackage[utf8]{inputenc}
\usepackage[english]{babel}
\usepackage{csquotes}

\usepackage{biblatex}
\renewbibmacro{in:}{%
        \ifentrytype{article}{}{\printtext{\bibstring{in}\intitlepunct}}}

\usepackage[hang,flushmargin]{footmisc} 

\usepackage{amsmath}
\usepackage{amssymb}
\usepackage{amsthm}

\usepackage{abstract}
\addto\captionsenglish{}

\newtheorem{defn}{Definition}[section]
\newtheorem{theo}[defn]{Theorem}
\newtheorem{prop}[defn]{Proposition}
\newtheorem{lemma}[defn]{Lemma}
\newtheorem{cor}[defn]{Corollary}

\newcommand{\Z}{\mathbb{Z}}

\newcommand{\R}{\mathbb{R}}

\newcommand{\func}[3]{#1:#2\rightarrow#3}
\newcommand{\norm}[2]{{||#1||}_{#2}}

\renewcommand{\epsilon}{\varepsilon}
\renewcommand{\theta}{\vartheta}
\renewcommand{\phi}{\varphi}

\title{Compactivorous Sets in Banach Spaces}

\author{Davide Ravasini}
\address{Institut für Mathematik, Universität Innsbruck \newline\indent Technikerstraße 13, 6020 Innsbruck, Austria}
\email{davide.ravasini@uibk.ac.at}

\begin{document}
\maketitle 
\let\thefootnote\relax\footnote{\today \newline \indent\emph{2020 Mathematics Subject Classification:} 46B20 (primary), 46B50, 54H11 (secondary). \newline
}

\begin{abstract}
\noindent \textsc{Abstract}. A set $E$ in a Banach space $X$ is compactivorous if for every compact set $K$ in $X$ there is a nonempty, (relatively) open subset of $K$ which can be translated into $E$. In a separable Banach space, this is a sufficient condition which guarantees the Haar nonnegligibility of Borel subsets.  We give some characterisations of this property in both separable and nonseparable Banach spaces and prove an extension of the main theorem to countable products of locally compact Polish groups.
\end{abstract}

\section{Introduction}
The purpose of the present paper is to study the notion of compactivorous sets in the context of Banach spaces. Let $E$ be a set in a Banach space $X$. $E$ is said to be \emph{compactivorous} if, for any compact set $K\subset X$, there exist an open set $V$ in $X$ and $x\in X$ such that $K\cap V\ne\varnothing$ and $x+K\cap V\subset E$. The first mention of this property can be found in \cite{emm_2016}, where it is provided as a sufficient condition which ensures that a given Borel set in a separable Banach space is not Haar null. A Borel set $E$ in a separable Banach space $X$ is \emph{Haar null} if there exists a Borel probability measure $\mu$ with compact support on $X$ such that $\mu(x+E)=0$ for all $x\in X$. It is almost immediate to see why Borel compactivorous sets cannot be Haar null. Let $E$ be compactivorous and consider a Borel probability measure $\mu$ with compact support $K\subset X$. By definition, there is an open set $V$ and $x\in X$ such that $K\cap V\ne\varnothing$ and $x+K\cap V\subset E$. Since $K\cap V$ is open and nonemtpy in $K$, $\mu(K\cap V)>0$, which implies $\mu(-x+E)>0$.

Haar null sets were first defined in every Abelian Polish group by Christensen in \cite{chris_1972} and research on them has been going on for decades. The survey articles \cite{bog_2018}, \cite{en_2020}, as well as \cite{benlin}, Chapter 6, provide a thorough introduction to the topic. Our study of compactivorous sets is motivated by the following question: let $C$ be a closed, convex set in a separable Banach space $X$. If $C$ is not Haar null, must it be compactivorous? To the best of the author's knowledge, the question is still unanswered. It is also worth mentioning that Darji introduced in \cite{dar_2013} a categorical version of Haar null sets, which for separable Banach spaces reads as follows: a Borel set $E$ in a separable Banach space $X$ is Haar meagre if there exists a compact metric space $M$ and a continuous function $\func{f}{M}{X}$ such that $f^{-1}(x+E)$ is meagre in $M$ for all $x\in X$. Observe that closed subsets of $X$ are compactivorous if and only if they are not Haar meagre. Our aim is to show some characterisations of compactivorous sets in (not necessarily separable) Banach spaces. In particular, we prove that, given a compactivorous set $E$, there is $r>0$ such that $E$ contains a translation of every compact set which lies in a closed ball of radius $r$, leading therefore to a seemingly stronger property.

Sections \ref{sec:zoom} and \ref{sec:sat} are devoted to the presentation of two of the main tools that will be needed in the proof of the main theorem, which is developed in Section \ref{sec:main}. Section \ref{sec:groups} shows a way to extend the main theorem to countable products of arbitrary, locally compact Polish groups. The standard notation of Banach space theory is used. If $X$ is a Banach space, $B_X$ denotes the closed unit ball of $X$, whereas $B(x,\epsilon)$ denotes the open ball centered in $x$ and with radius $\epsilon$. Given a set $S$, $\overline{\textup{span}}\,S$ denotes the smallest closed subspace which contains it. The interior and closure of a set $S$ in a topological space are denoted $\textup{int}\,S$ and $\textup{cl}\,S$ respectively and $2^S$ stands for the power set of a set $S$. $\ell_1$ is the Banach space of absolutely summable sequences, whereas $c_0$ is the Banach space of sequences which converge to $0$. All Banach spaces are assumed to be real.

\section{The rescaling property}
\label{sec:zoom}
The following, simple property of compactivorous sets works as a motivation for the definition we are going to introduce later.
\begin{lemma}
\label{lemma:resc1}
Let $E$ be a compactivorous set in a Banach space. Then $E$ is compactivorous by rescaling. That is: for each compact set $K\subset X$ there are $x\in X$ and $\delta>0$ such that $x+\delta K\subset E$.
\end{lemma}
\begin{proof}
Let $K$ be a compact subset of $X$. Since every compact set is contained in a convex compact set, we can assume that $K$ is convex. Pick an open set $V$ and $x\in X$ such that $K\cap V\ne\varnothing$ and $x+K\cap V\subset E$. Without loss of generality, we can assume that $0\in K\cap V$ and that $V=B(0,\epsilon)$ for some $\epsilon>0$. Let $\delta\in(0,1)$ be such that $\delta K\subset V$. Since $K$ is convex, $\delta K\subset K$, hence $\delta K\subset K\cap V$. Thus, $x+\delta K\subset E$.
\end{proof}

Let $\mathcal{F}$ be a family of subsets in a separable Banach space $X$. Drawing inspiration from Lemma \ref{lemma:resc1}, we say that $\mathcal{F}$ has the \emph{rescaling property} if the two following conditions are satisfied.
\begin{enumerate}
\item $\mathcal{F}$ is hereditary: if $A\in\mathcal{F}$ and $B\subseteq A$, then $B\in\mathcal{F}$.
\item For each compact set $K\subset X$, there are $\delta>0$ and $x\in X$ such that $x+\delta K\in\mathcal{F}$.
\end{enumerate}
\begin{lemma}
\label{lemma:resc2}
Let $E$ be a subset of a Banach space $X$. $E$ is compactivorous by rescaling if and only if every closed, separable subspace $Y\subseteq X$ admits a family of subsets with the rescaling property $\mathcal{F}_Y$ such that every $S\in\mathcal{F}_Y$ can be translated into $E$.
\end{lemma}
\begin{proof}
Suppose that $E$ is compactivorous by rescaling and let $Y$ be a separable subspace of $X$. For each compact set $K\subset Y$, choose $\delta_K>0$ such that $\delta_K\cdot K$ can be translated into $E$ and define
\[ \mathcal{F}_Y=\{S\subset Y\,:\,\text{there is a compact set }K\subset Y\text{ such that }S\subseteq\delta_K\cdot K\}. \]
$\mathcal{F}_Y$ is the family with the rescaling property we were looking for.

Conversely, let $K\subset X$ be compact and define $Y=\overline{\textup{span}}\,K$. $Y$ is a separable subspace of $X$ and admits therefore a family of subsets $\mathcal{F}_Y$ with the rescaling property such that every $S\in\mathcal{F}_Y$ can be translated into $E$. Hence, there are $y\in Y$ and $\delta>0$ such that $y+\delta K\in\mathcal{F}_Y$ and $z\in X$ such that $z+y+\delta K\subset E$. Setting $x=z+y$, we have $x+\delta K\subset E$ and since $K$ is arbitrary, this concludes the proof.
\end{proof}

As we shall see later, Lemma \ref{lemma:resc2} allows to reduce the study of compactivorous sets to separable Banach spaces by looking at families with the rescaling property instead of compactivorous sets themselves.

\section{Saturated compact subsets of $\ell_1$}
\label{sec:sat}
The second tool we need is the concept of saturated compact subsets of $\ell_1$. To have a first understanding of what they are it is important to recall the following, well-known characterisation of compact subsets of $\ell_1$.
\begin{prop}
\label{prop:comp1}
A set $S\subset\ell_1$ is totally bounded if and only if
\begin{enumerate}
\item $S$ is pointwise bounded, i.e. for any $n\geq 1$ the set $\{ x(n):x\in S\}$ is bounded;
\item For any $\epsilon>0$ there is $n_0\geq 1$ such that
\[ \sum_{j=n}^{\infty}|x(j)|<\epsilon \]
for any $x\in S$ and $n\geq n_0$. 
\end{enumerate}
In particular, a set $K\subset\ell_1$ is compact if and only if $K$ is closed and satisfies conditions (1) and (2).
\end{prop}

Proposition \ref{prop:comp1} suggests we introduce the following definition: given $x\in\ell_1$, the \emph{modulus of decay} of $x$ is the sequence $\omega_x\in c_0$ given by
\[ \omega_x(n)=\sum_{j=n}^{\infty}|x(j)|. \]
Notice that, for each $x\in\ell_1$, the sequence $\omega_x$ is nonnegative and decreasing, that $\omega_x(1)=\norm{x}{\ell_1}$ and that the function $\func{\phi}{\ell_1}{c_0}$ which assigns to each $x\in\ell_1$ its modulus of decay is continuous. If $K\subset\ell_1$ is compact, Proposition \ref{prop:comp1} guarantees that the sequence $\Omega_K$ given by
\[ \Omega_K(n)=\sup_{x\in K}\omega_x(n) \]
belongs to $c_0$, and we shall call it the \emph{modulus of decay} of the set $K$. We say that a compact set $K\subset\ell_1$ is \emph{saturated} if it consists exactly of all elements $x\in\ell_1$ such that $\omega_x\leq\Omega_K$. Saturated compact sets are clearly symmetric and convex. However, they are far from being a countable product of intervals. Consider for instance
\[ K=\prod_{n=1}^{\infty}\biggl[-\frac{1}{n^2},\frac{1}{n^2}\biggr]\subset\ell_1. \]
$K$ is a compact set whose modulus of decay is 
\[ \Omega_K(1)=\frac{\pi^2}{6},\quad\Omega_K(n)=\frac{\pi^2}{6}-\sum_{j=1}^{n-1}\frac{1}{j^2}\;\text{ if }n>1, \]
but the element $x=(\pi^2/6,0,0,\dots)$ does not belong to $K$, even though $\omega_x\leq\Omega_K$.

Finally, it is also important to keep in mind that to every nonnegative and decreasing sequence $z\in c_0$ there corresponds a saturated compact subset of $\ell_1$ whose modulus of decay is $z$, namely the set $K(z)=\{x\in\ell_1\,:\,\omega_x\leq z\}$. $K(z)$ is compact by Proposition \ref{prop:comp1} and $\Omega_{K(z)}=z$ by definition. In particular, this implies that every compact set is contained in a saturated compact set with the same modulus of decay.

To conclude this preliminary part, the next lemma describes what happens if we intersect a saturated compact set with a closed ball centered around the origin.
\begin{lemma}
\label{lemma:comp2}
Let $K\subset\ell_1$ be a saturated compact set with modulus of decay $\Omega_K$. For any $r>0$, the compact set $L=rB_{\ell_1}\cap K$ is also saturated and its modulus of decay $\Omega_L$ is given exactly by $\Omega_L(n)=\min\{r,\Omega_K(n)\}$.
\end{lemma}
\begin{proof}
Let $x\in\ell_1$ be such that $\omega_x(n)\leq\min\{r,\Omega_K(n)\}$. Then $x\in K$ because $K$ is saturated, and $\norm{x}{\ell_1}=\omega_x(1)\leq r$. Thus $x\in rB_{\ell_1}\cap K$. Conversely, if $x\in L$, then $\omega_x(n)\leq\Omega_K(n)$ for every $n$ because $L\subseteq K$, moreover $\omega_x(n)\leq\omega_x(1)\leq r$, therefore $\omega_x(n)\leq\min\{r,\Omega_K(n)\}$ for each $n$.
\end{proof}

\section{The main theorem}
\label{sec:main}
We are now ready to work towards a characterisation of compactivorous sets. As anticipated, we are actually going to study families with the rescaling property and perform a translation of the obtained results back in the original setting of compactivorous sets. We start with a lemma.
\begin{lemma}
\label{lemma:corr}
Let $X$ be a Banach space, let $Y$ be a closed subspace and denote with $\func{\pi}{X}{X/Y}$ the quotient map. For each $r>0$, the correspondence
\[ \func{f_r}{\textup{int}(rB_{X/Y})}{2^X},\quad f_r(x)=\pi^{-1}(x)\cap rB_X, \]
assumes closed, convex and nonempty values and it is lower hemicontinuous.
\end{lemma}
\begin{proof}
The only nontrivial assertion to verify is that $f_r$ is lower hemicontinuous for every $r>0$. To see this, define a second correspondence 
\[ \func{g_r}{\textup{int}(rB_{X/Y})}{2^X},\quad g_r(x)=\pi^{-1}(x)\cap\textup{int}(rB_X). \]
$g_r$ has nonempty images because $\textup{int}(rB_{X/Y})=\pi\bigl(\textup{int}(rB_X)\bigr)$ and, since the correspondence $\pi^{-1}$ is lower hemicontinuous and $\textup{int}(rB_X)$ is open, $g_r$ is also lower hemicontinuous. We want to show that $f_r$ is the closure of $g_r$, i.e.\ that $\textup{cl}\bigl(g_r(x)\bigr)=f_r(x)$ for each $x\in\textup{int}(rB_{X/Y})$. Then the lower hemicontinuity of $f$ follows from \cite{hitchhiker}, Lemma 17.22. Clearly, we have $\textup{cl}\bigl(\pi^{-1}(x)\cap\textup{int}(rB_X)\bigr)\subseteq\pi^{-1}(x)\cap rB_X$ for every $x\in\textup{int}(rB_{X/Y})$. Let $y\in\pi^{-1}(x)\cap rB_X$ and fix a vector $z\in\pi^{-1}(x)\cap\textup{int}(rB_X)$. For all positive integers $n$, define the elements $x_n=(1-2^{-n})y+2^{-n}z$. It is easy to check that $x_n\in\pi^{-1}(x)\cap\textup{int}(rB_X)$ for each $n$ and the sequence $(x_n)$ converges to $y$, hence $y\in\textup{cl}\bigl(\pi^{-1}(x)\cap\textup{int}(rB_X)\bigr)$ and the equality $\textup{cl}\bigl(g(x)\bigr)=f(x)$ holds.
\end{proof}

The next theorem characterises families with the rescaling property in a very strong sense.
\begin{theo}
\label{theo:main1}
Let $X$ be a separable Banach space and let $\mathcal{F}$ be a hereditary family of subsets of $X$. The two following assertions are equivalent.
\begin{enumerate} 
\item $\mathcal{F}$ has the rescaling property.
\item There is $r>0$ with the property that, for every compact set $K\subset rB_X$, one can find an $x\in X$ such that $x+K\in\mathcal{F}$.
\end{enumerate}
\end{theo}
\begin{proof}
First, let us observe that (2)$\Rightarrow$(1) is obvious, because for every compact set $K$ there is $\delta>0$ such that $\delta K$ fits inside $rB_X$. Further, let $\func{T}{X}{Z}$ be an isomorphism between $X$ and another separable Banach space $Z$ and consider the hereditary family of subsets $T(\mathcal{F})=\{T(S)\,:\,S\in\mathcal{F}\}$. It is trivial to check that $\mathcal{F}$ satisfies condition (1) in the statement of the theorem if and only if $T(\mathcal{F})$ does, and the same holds for condition (2). Thus, it suffices to focus our attention on quotient spaces of $\ell_1$. From now on, we assume that $X=\ell_1/Y$, where $Y$ is a closed subspace of $\ell_1$, and that $\func{\pi}{\ell_1}{X}$ is the quotient map. 

Arguing by way of contradiction, suppose that, for each positive integer $n$, there is a compact set $K_n\subset \textup{int}(n^{-1}B_X)$ such that no translation of $K_n$ belongs to $\mathcal{F}$. Consider for each $n$ the correspondence
\[ \func{f_n}{\textup{int}(n^{-1}B_{X})}{2^{\ell_1}},\quad f_n(x)=\pi^{-1}(x)\cap n^{-1}B_{\ell_1}, \]
which is lower hemicontinuous by Lemma \ref{lemma:corr}. For each $n$, let $s_n$ be a continuous selection of $f_n$, which exists thanks to Michael's selection theorem (see \cite{mich_1956}, Theorem 3.2$''$) and observe that $s_n(K_n)$ is a compact set in $n^{-1}B_{\ell_1}$. Define
\[ K=\{0\}\cup\bigcup_{n=1}^{\infty}s_n(K_n). \]
$K$ is compact and it comes therefore with its modulus of decay $\Omega_K$. Let $L$ be the saturated compact set with the same modulus of decay and let $H\subset\ell_1$ be another saturated compact set whose modulus of decay $\Omega_H$ satisfies
\[ \lim_{j\to\infty}\frac{\Omega_H(j)}{\Omega_K(j)}=+\infty. \]
Choose $\delta>0$ and observe that $\delta H$ is another saturated compact set whose modulus of decay is $\delta\Omega_H$. By our choice of $H$, there is the smallest index $j_0$ such that $\delta\Omega_H(j)>\Omega_K(j)$ for all $j\geq j_0$. If we set $t=\Omega_K(j_0)$ and let $y\in tB_{\ell_1}\cap L$, then by Lemma \ref{lemma:comp2} we have $\omega_y(j)\leq\min\{t,\Omega_K(j)\}\leq\delta\Omega_H(j)$ for all positive integers $j$, which means that $tB_{\ell_1}\cap L\subseteq\delta H$. Let $n$ be such that $n^{-1}\leq t$ and observe that $s_n(K_n)\subseteq n^{-1}B_{\ell_1}\cap L\subseteq tB_{\ell_1}\cap L\subseteq\delta H$. If we apply the quotient map, we get that $K_n=\pi\bigl(s_n(K_n)\bigr)\subseteq\delta\pi(H)$. This implies that there is no $x\in X$ such that $x+\delta\pi(H)\in\mathcal{F}$, otherwise we would conclude that $x+K_n\in\mathcal{F}$, because $\mathcal{F}$ is hereditary. As $\pi(H)$ is compact and $\delta$ is arbitrary, this shows that (1) fails.
\end{proof}

We are finally ready to combine all the results we have obtained so far in the characterisation theorem for compactivorous sets.
\begin{theo} 
\label{theo:main2}
Let $E$ be a subset in a Banach space $X$. Then the following assertions are equivalent. 
\begin{enumerate}
\item $E$ is compactivorous.
\item For every compact set $K\subset X$, there are $x\in X$ and $\delta>0$ such that $x+\delta K\subset E$.
\item There is $r>0$ with the property that, for every compact set $K\subset rB_X$, one can find an $x\in X$ such that $x+K\subset E$.
\end{enumerate}
\end{theo}
\begin{proof}
The implication (1)$\Rightarrow$(2) is the content of Lemma \ref{lemma:resc1}, whereas (3)$\Rightarrow$(1) is obvious: if condition (3) is satisfied, then it suffices to intersect every compact set with a sufficiently small open ball to prove that $E$ is compactivorous. It remains to show that (2)$\Rightarrow$(3). 

If $E$ satisfies (2) and $Y$ is a closed, separable subspace of $X$, then by Lemma \ref{lemma:resc2} there is a family $\mathcal{F}$ of subsets of $Y$ with the rescaling property such that every element of $\mathcal{F}$ can be translated into $E$. By Theorem \ref{theo:main1}, there is $r>0$ with the property that, for each compact set $K\subset rB_Y$, there is $y\in Y$ such that $y+K\in\mathcal{F}$. In particular, every compact subset of $rB_Y$ can be translated into $E$. This justifies the introduction of the following notation: given a closed, separable subspace $Y\subseteq X$, we denote by $\rho(E,Y)$ the supremum of all $r>0$ such that every compact subset of $rB_Y$ can be translated into $E$. Observe that if $Y$ and $Z$ are closed, separable subspaces of $X$ such that $Y\subseteq Z$, then $\rho(E,Z)\leq\rho(E,Y)$. Since every compact set lies in a closed, separable subspace, to ensure that $E$ satisfies (3) it remains to prove that there is $r>0$ such that $\rho(E,Y)>r$ for all closed, separable subspaces $Y$. Arguing by contradiction, suppose that for each positive integer $n$ there is a closed, separable subspace $Y_n$ such that $\rho(E,Y_n)\leq n^{-1}$. Define
\[ Y=\overline{\textup{span}}\,\biggl(\,\bigcup_{n=1}^\infty Y_n\biggr). \]
$Y$ is a closed, separable subspace of $X$ and therefore it has its $\rho(E,Y)>0$. At the same time though, we have $Y\supseteq Y_n$ for each $n$, hence $\rho(E,Y)\leq\rho(E,Y_n)\leq n^{-1}$ for each $n$, which is nonsense.
\end{proof}

\begin{cor}
Let $X$ be a Banach space, $Y$ a closed subspace and $E$ a subset of $X/Y$. Let $\func{\pi}{X}{X/Y}$ be the quotient map. Then $E$ is compactivorous if and only if $\pi^{-1}(E)$ is compactivorous.
\end{cor}
\begin{proof}
Assume that $\pi^{-1}(E)$ is compactivorous and let $s$ be a continuous selection of the correspondence $\func{\pi^{-1}}{X/Y}{2^X}$ (we use once again \cite{mich_1956}, Theorem 3.2$''$). If $K\subset X/Y$ is compact, so is $s(K)$, hence there are $x\in X$ and $\delta>0$ such that $x+\delta s(K)\subset\pi^{-1}(E)$. Applying $\pi$, one gets $\pi(x)+\delta K\subset E$. As $K$ is arbitrary, this shows that $E$ is compactivorous. Conversely, let $E$ be compactivorous and let $K\subset X$ be compact. Then there are $x\in X$ and $\delta>0$ such that $\pi(x)+\delta\pi(K)\subset E$, which implies $x+\delta K\subset\pi^{-1}(E)$. The conclusion follows once again because $K$ is arbitrary.
\end{proof}

\section{Generalising to topological groups}
\label{sec:groups}
The notion of a compactivorous set can be extended to arbitrary Hausdorff topological groups. Let $G$ be a Hausdorff topological group, let $e$ be its unit and let $E$ be a subset of $G$. We say that $E$ is compactivorous if for every compact set $K\subseteq G$ there are an open set $V\subseteq G$ and $g,h\in G$ such that $K\cap V\ne\varnothing$ and $g(K\cap V)h\subseteq E$. We shall also say that $E$ is \emph{strongly compactivorous} if there is a closed neighbourhood $U$ of $e$ such that, for every compact set $K\subseteq U$, there are $g,h\in G$ such that $gKh\subseteq E$. Clearly, every strongly compactivorous set is compactivorous. Let us say that $G$ is \emph{fattening} if also the converse holds, i.e.\ if every compactivorous subset of $G$ is strongly compactivorous. Locally compact groups are trivially fattening and that every Banach space is fattening is the content of Theorem \ref{theo:main2}. The next example shows that there are in fact nonfattening groups.

\begin{prop}
Let $X$ be an infinite-dimensional Banach space with the Schur property. Then $X$ is nonfattening when endowed with the weak topology.
\end{prop}
\begin{proof}
We show that the unit ball $B_X$ is compactivorous but not strongly compactivorous. Let $K$ be a weakly compact set in $X$. Since $X$ has the Schur property, by the Eberlein-\v{S}mulian theorem $K$ is also norm-compact. The restrictions of both the weak and the norm topology to $K$ make it a compact, Hausdorff space, hence they must agree. $B_X$ has nonempty interior in the norm topology, which means that it is compactivorous in the norm topology. Let $U$ be a norm-open set such that $K\cap U$ is nonempty and can be translated into $B_X$. Let $V$ be a weakly open set such that $K\cap U=K\cap V$. Then $K\cap V$ can be translated into $B_X$ and, since $K$ is arbitrary, this shows that $B_X$ is compactivorous in the weak topology. However, $B_X$ cannot be strongly compactivorous, as every neighbourhood of $0$ is unbounded and contains therefore compact sets of arbitrarily large diameter. 
\end{proof}

Finally, we prove that countable products of locally compact Polish groups are fattening.

\begin{theo}
Let ${\{G_n\}}_{n=1}^{\infty}$ be a countable family of locally compact Polish groups and set
\[ G=\prod_{n=1}^{\infty}G_n. \]
Then $G$ is fattening.
\end{theo}
\begin{proof}
We start by fixing some notation. Let $d$ be a complete metric which generates the topology of $G$ and, for each $n\geq 1$, let $\func{p_n}{G}{G_n}$ be the usual projection map. Let $e_n$ be the unit of $G_n$ and $e=(e_1,e_2,\dots)$ be the unit of $G$. Let $E$ be a compactivorous subset of $G$, set $B_k=\{g\in G\,:\,d(g,e)\leq k^{-1}\}$ for each positive integer $k$ and, looking for a contradiction, assume that for each $k$ there is a compact set $K_k\subseteq B_k$ such that $aK_kb$ does not lie entirely in $E$ for any $a,b\in G$. The set
\[ K=\{e\}\cup\bigcup_{k=1}^\infty K_k \]
is compact and so are the sets $p_n(K)$ for all $n\geq 1$. Since $G_n$ is locally compact, there is an open neighbourhood $U_n$ of $p_n(K)$ such that $\textup{cl}(U_n)$ is compact. Define the compact set
\[ L=\prod_{n=1}^{\infty}\textup{cl}(U_n). \]
Since $E$ is compactivorous, there must be a positive integer $m$, $h_1,h_2\in G$ and an open set $V\subseteq G$ of the form $V=V_1\times\cdots\times V_m\times G_{m+1}\times\cdots$, where $V_n\subseteq G_n$ is open for all $n\in\{1,\dots,m\}$, such that $L\cap V\ne\varnothing$ and $h_1(L\cap V)h_2\subseteq E$. In particular, $U_n\cap V_n\ne\varnothing$ for all $n\in\{1,\dots,m\}$. For every $n\in\{1,\dots,m\}$, choose $g_n\in U_n\cap V_n$ and observe that $g_n^{-1}(U_n\cap V_n)$ is an open neighbourhood of $e_n$. Moreover, ${\{p_n(B_k)\}}_{k=1}^{\infty}$ is a nested neighbourhood basis of $e_n$ for every $n\geq 1$, hence for each $n\in\{1,\dots,m\}$ there exists $k_0(n)$ such that $p_n(B_k)\subseteq g_n^{-1}(U_n\cap V_n)$ for all $k\geq k_0(n)$, which implies that $p_n(K_k)\subseteq g_n^{-1}(U_n\cap V_n)$ for all $k\geq k_0(n)$. Let $j$ be a sufficiently large index such that $p_n(K_j)\subseteq g_n^{-1}(U_n\cap V_n)$ for all $n\in\{1,\dots,m\}$ and set $g=(g_1,\dots,g_m,e_{m+1},\dots)$. Then
\begin{align*} 
gK_j & \subseteq g\biggl(\,\prod_{n=1}^{\infty}p_n(K_j)\biggr)=\prod_{n=1}^{m}g_np_n(K_j)\times\prod_{n=m+1}^\infty p_n(K_j)\subseteq \\
& \subseteq\prod_{n=1}^m(U_n\cap V_n)\times\prod_{n=m+1}^{\infty}\textup{cl}(U_n)\subseteq L\cap V.
\end{align*}
However, this is a contradiction, as it implies that $h_1gK_jh_2\subseteq E$, against the assumption on $K_j$.
\end{proof}

\begin{cor}
The Fréchet space $\R^\omega$ and the Baire space $\Z^\omega$ are fattening.
\end{cor}

\section{Open questions}
The following questions arose naturally during the writing of this paper.
\begin{enumerate}
\item What other Hausdorff topological groups are fattening? What can be said about the weak topologies of Banach spaces without the Schur property or about Polish groups in general? 
\item Esterle, Matheron and Moreau were the first to ask in \cite{emm_2016} whether every closed, convex set which is not Haar null in a separable Banach space must be compactivorous. We observe that it is enough to ask the same question in $\ell_1$. Indeed, let $Y$ be a closed subspace of $\ell_1$, let $C\subset \ell_1/Y$ be a closed, convex subset and let $\func{\pi}{\ell_1}{\ell_1/Y}$ be the quotient map. Provided that the answer is affirmative, one obtains the following chain of implications: $C$ is not Haar null $\Rightarrow\pi^{-1}(C)$ is not Haar null $\Rightarrow\pi^{-1}(C)$ is compactivorous $\Rightarrow C$ is compactivorous. Using invariance under isomorphisms, the result would then become true in all separable Banach spaces.
\end{enumerate}

\section*{Acknowledgements}
The author would like to thank professors Eva Kopeck\'{a} and Christian Bargetz for their precious suggestions. The work of the author is supported by the Austrian Science Fund (FWF): P 32523-N.

\printbibliography
\end{document}

This section is devoted to presenting some of the machinery that is involved in the proof of the main result. 

Let $X$ be a separable Banach space. A family $\mathcal{F}$ of subsets of $X$ is said to have the \emph{compactivorous property} or to be compactivorous if the following two conditions are satisfied.
\begin{enumerate}
\item $\mathcal{F}$ is hereditary: if $A\in\mathcal{F}$ and $B\subseteq A$, then $B\in\mathcal{F}$.
\item For every compact set $K\subset X$, there are an open set $V\subseteq X$ and $x\in X$ such that $K\cap V\ne\varnothing$ and $x+K\cap V\in\mathcal{F}$.
\end{enumerate}
The fact that a given set in a Banach space is compactivorous can be expressed by using families with the compactivorous property.
\begin{prop}
\label{prop:fam}
Let $X$ be a Banach space. A set $E\subset X$ is compactivorous if and only if every closed, separable subspace $Y\subseteq X$ admits a family of subsets $\mathcal{F}_Y$ with the compactivorous property such that every $S\in\mathcal{F}_Y$ can be translated into $E$.
\end{prop}
\begin{proof}
Suppose that $E$ is compactivorous and let $Y$ be a separable subspace of $X$. For each compact set $K\subset Y$ choose an open set $U_K\subseteq X$ such that $K\cap U_K\ne\varnothing$ and $K\cap U_K$ can be translated into $E$. Define $V_K=U_K\cap Y$ for each $K$ and set
\[ \mathcal{F}_Y=\{S\subset Y\,:\,\text{there is a compact set }K\subset Y\text{ such that }S\subseteq K\cap V_K\}. \]
$\mathcal{F}_Y$ is a family with the compactivorous property consisting of sets which can be translated into $E$.

Conversely, assume that every separable subspace admits such a family and let $K\subset X$ be compact. Define $Y=\overline{\textup{span}}\,K$ and observe that $Y$ is separable. Let $\mathcal{F}_Y$ be a family of subsets of $Y$ with the compactivorous property made of elements which can be translated into $E$. Let $V$ be an open set in $Y$ and let $y\in Y$ be such that $K\cap V\ne\varnothing$ and $y+K\cap V\in\mathcal{F}_Y$. Consider an open set $U$ in $X$ such that $U\cap Y=V$ and pick $x\in X$ such that $x+y+K\cap V\subset E$. If we set $z=x+y$, then $z+K\cap U\subset E$. Since $K$ is arbitrary, $E$ is compactivorous.
\end{proof}

\noindent As we are going to see later, proposition \ref{prop:fam} allows to reduce the problem of studying compactivorous sets to separable spaces by looking at families with the compactivorous property instead of compactivorous sets themselves.

\section{Characterisation of compactivorous families}
We are now ready to work towards a characterisation of compactivorous sets. As anticipated, we are actually going to study families with the rescaling property and perform a translation of the obtained results back in the original setting of compactivorous sets. We start with a lemma.
\begin{lemma}
\label{lemma:corr}
Let $X$ be a Banach space, let $Y$ be a closed subspace and denote with $\func{\pi}{X}{X/Y}$ the quotient map. For each $r>0$, the correspondence
\[ \func{f_r}{\textup{int}(rB_{X/Y})}{2^X},\quad f_r(x)=\pi^{-1}(x)\cap rB_X, \]
assumes closed, convex and nonempty values and it is lower hemicontinuous.
\end{lemma}
\begin{proof}
The only nontrivial assertion to verify is that $f_r$ is lower hemicontinuous for every $r>0$. To see this, define a second correspondence 
\[ \func{g_r}{\textup{int}(rB_{X/Y})}{2^X},\quad g_r(x)=\pi^{-1}(x)\cap\textup{int}(rB_X). \]
$g_r$ has nonempty images because $\textup{int}(rB_{X/Y})=\pi\bigl(\textup{int}(rB_X)\bigr)$ and, since the correspondence $\pi^{-1}$ is lower hemicontinuous and $\textup{int}(rB_X)$ is open, $g_r$ is also lower hemicontinuous. We want to show that $f_r$ is the closure of $g_r$, i.e.\ that $\textup{cl}\bigl(g_r(x)\bigr)=f_r(x)$ for each $x\in\textup{int}(rB_{X/Y})$. Then the lower hemicontinuity of $f$ follows from \cite{hitchhiker}, Lemma 17.22. Clearly, we have $\textup{cl}\bigl(\pi^{-1}(x)\cap\textup{int}(rB_X)\bigr)\subseteq\pi^{-1}(x)\cap rB_X$ for every $x\in\textup{int}(rB_{X/Y})$. Let $y\in\pi^{-1}(x)\cap rB_X$ and fix a vector $z\in\pi^{-1}(x)\cap\textup{int}(rB_X)$. For all positive integers $n$, define the elements $x_n=(1-2^{-n})y+2^{-n}z$. It is easy to check that $x_n\in\pi^{-1}(x)\cap\textup{int}(rB_X)$ for each $n$ and the sequence $(x_n)$ converges to $y$, hence $y\in\textup{cl}\bigl(\pi^{-1}(x)\cap\textup{int}(rB_X)\bigr)$ and the equality $\textup{cl}\bigl(g(x)\bigr)=f(x)$ holds.
\end{proof}

\begin{theo}
\label{theo:famcomp}
Let $X$ be a separable Banach space and let $\mathcal{F}$ be a hereditary family of subsets of $X$. The following conditions are equivalent.
\begin{enumerate}
\item $\mathcal{F}$ has the compactivorous property;
\item There is $r>0$ with the property that, for each compact set $K\subset rB_X$, one can find $x\in X$ such that $x+K\in\mathcal{F}$;
\item For every compact set $K\subset X$ there are $\delta>0$ and $x\in X$ such that $x+\delta K\in\mathcal{F}$.
\end{enumerate}
\end{theo}
\begin{proof}
First, let us observe that the implication (2)$\Rightarrow$(3) is obvious. We shall show that (1)$\Rightarrow$(2) and that (3)$\Rightarrow$(1). Further, let $\func{T}{X}{Z}$ be an isomorphism between $X$ and another separable Banach space $Z$ and consider the hereditary family of subsets $T(\mathcal{F})=\{T(S)\,:\,S\in\mathcal{F}\}$. For $i\in\{1,2,3\}$, it is trivial to check that $\mathcal{F}$ satisfy condition $(i)$ in the statement of the theorem if and only if $T(\mathcal{F})$ does, therefore it suffices to focus our attention on quotient spaces of $\ell_1$. From now on, we assume that $X=\ell_1/Y$, where $Y$ is a closed subspace of $\ell_1$, and that $\func{\pi}{\ell_1}{X}$ is the quotient map.

(1)$\Rightarrow$(2) Arguing by way of contradiction, suppose that, for each positive integer $n$, there is $K_n\subset \textup{int}(n^{-1}B_X)$ such that no translation of $K_n$ belongs to $\mathcal{F}$. Consider for each $n$ the correspondence
\[ \func{f_n}{\textup{int}(n^{-1}B_{X})}{2^{\ell_1}},\quad f_n(x)=\pi^{-1}(x)\cap n^{-1}B_{\ell_1}, \]
which is lower hemicontinuous by Lemma \ref{lemma:corr}. Let $s_n$ be a continuous selection of $f_n$, which exists thanks to Michael's selection theorem (see \cite{mich_1956}, Theorem 3.2'') and observe that $s_n(K_n)$ is a compact set in $n^{-1}B_{\ell_1}$. Define
\[ K=\{0\}\cup\bigcup_{n=1}^{\infty}s_n(K_n). \]
$K$ is compact and it comes therefore its modulus of decay $\Omega_K$. Let $L$ be the saturated compact set with the same modulus of decay. Let $V$ be any open set such that $\pi(L)\cap V\ne\varnothing$, which implies that $U=\pi^{-1}(V)$ is an open set and that $L\cap U\ne\varnothing$. Let $z\in c_{00}\cap L\cap U$ (we can pick such $z$ because $c_{00}\cap L$ is dense in $L$) and let $k_0$ be the minimal index such that $z(k)=0$ for all $k>k_0$. Since $z$ lies in the interior of $U$, we can slightly decrease the absolute value of $z(k_0)$ if necessary in order to safely assume that $\omega_z(k)<\Omega_K(k)$ for all $k\in\{1,\dots,k_0\}$. Let $\epsilon>0$ be such that $B(z,\epsilon)\subseteq U$. Set $\delta=\min\{\Omega_K(k)-\omega_z(k):1\leq k\leq k_0\}$ and choose $n\geq 1$ such that $n^{-1}<\min\{\delta,\epsilon\}$. Let $y\in s_n(K_n)$. Since $s_n(K_n)\subseteq n^{-1}B_{\ell_1}\cap L$, by Lemma \ref{lemma:comp2} one has $\omega_y(k)\leq n^{-1}$ for all $k\in\{1,\dots,k_0\}$, therefore
\[ \omega_{y+z}(k)\leq\omega_y(k)+\omega_z(k)\leq n^{-1}+\omega_z(k)<\delta+\omega_z(k)\leq\Omega_K(k) \]
for all $k\in\{1,\dots,k_0\}$. If $k>k_0$, then $\omega_z(k)=0$ and it follows immediately that $\omega_{y+z}(k)\leq\Omega_K(k)$. This proves that $y+z\in L$. Finally, it is clear that $y+z\in B(z,\epsilon)$ for all $y\in s_n(K_n)$, thus $z+s_n(K_n)\subseteq L\cap U$. By applying the quotient map, we get $\pi(z)+K_n\subseteq\pi(L)\cap V$. This implies that there is no $x\in X$ such that $x+\pi(L)\cap V\in\mathcal{F}\,$: if there were such $x$, we would conclude that $x+\pi(z)+K_n\in\mathcal{F}$, because $\mathcal{F}$ is hereditary, against the assumption on $K_n$. As $V$ has been chosen arbitrarily, we have just proved that no translation of any open subset of $\pi(L)$ belongs to $\mathcal{F}$, hence (1) fails.

(3)$\Rightarrow$(1) 
\end{proof}

\section{Families with the compactivorous property}
Let $X$ be a separable Banach space. A family $\mathcal{F}$ of subsets of $X$ is said to have the \emph{compactivorous property} if it satisfies the following two conditions.
\begin{enumerate}
\item $\mathcal{F}$ is hereditary: if $A\in\mathcal{F}$ and $B\subseteq A$, then $B\in\mathcal{F}$.
\item For every compact set $K\subset X$, there are an open set $V\subseteq X$ and $x\in X$ such that $K\cap V\ne\varnothing$ and $x+K\cap V\in\mathcal{F}$.
\end{enumerate}
\begin{prop}
Let $X$ be a Banach space. A set $E\subset X$ is compactivorous if and only if every closed separable subspace $Y\subseteq X$ admits a family of subsets $\mathcal{F}_Y$ with the compactivorous property such that every $S\in\mathcal{F}_Y$ can be translated into $E$.
\end{prop}
\begin{proof}
Suppose that $E$ is compactivorous and let $Y$ be a separable subspace of $X$. For each compact set $K\subset Y$ choose an open set $U_K\subseteq X$ such that $K\cap U_K\ne\varnothing$ and $K\cap U_K$ can be translated into $E$. Define $V_K=U_K\cap Y$ for each $K$ and set
\[ \mathcal{F}_Y=\{S\subset Y\,:\,\text{there is a compact set }K\subset Y\text{ such that }S\subseteq K\cap V_K\}. \]
$\mathcal{F}_Y$ is a family with the compactivorous property consisting of sets which can be translated into $E$.

Conversely, assume that every separable subspace admits such a family and let $K\subset X$ be compact. Define $Y=\overline{\textup{span}}K$ and observe that $Y$ is separable. Let $\mathcal{F}_Y$ be a family of subsets of $Y$ with the compactivorous property made of elements which can be translated into $E$ and let $V$ be an open set in $Y$ such that $K\cap V\ne\varnothing$ and $y+K\cap V\in\mathcal{F}_Y$ for some $y\in Y$. Let $U$ be an open set in $X$ such that $U\cap X=V$ and let $x\in X$ be such that $x+y+K\cap V\subset E$. Then $x+y+K\cap U\subset E$. Since $K$ was arbitrary, $E$ is compactivorous.
\end{proof}

\begin{lemma}
Let $X$ be a Banach space and let ${\{Y_k\}}_{k=1}^\infty$ be a countable family of separable, closed subspaces. Then 
\[ Y=\overline{\textup{span}}\biggl(\,\bigcup_{k=1}^\infty Y_n\biggr) \] 
is separable.
\end{lemma}
\begin{proof}
For each $k$, let $D_k$ be a countable, dense subset of $Y_k$, define
\[ D=\biggl\{\,\sum_{k=1}^nz_k\,:\,n\geq1,\, z_k\in D_k\text{ for all }k\in\{1,\dots,n\}\biggr\} \]
and observe that $D$ is countable. Choose $y\in Y$ and fix $\epsilon>0$. Find $n\geq 1$ and $y_k\in Y_k$ for $k\in\{1,\dots,n\}$ such that 
\[ {\biggl\lVert y-\sum_{k=1}^ny_k\biggr\rVert}_X<\frac{\epsilon}{2}. \]
Finally, for each $k\in\{1,\dots,n\}$ find $z_k\in D_k$ such that ${\lVert y_k-z_k\rVert}_X<{(2n)}^{-1}\epsilon$. Then we have
\[ {\biggl\lVert y-\sum_{k=1}^nz_k\biggr\rVert}_X\leq{\biggl\lVert y-\sum_{k=1}^ny_k\biggr\rVert}_X+\sum_{k=1}^n{\lVert y_k-z_k\rVert}_X<\frac{\epsilon}{2}+n\cdot\frac{\epsilon}{2n}=\epsilon. \]
Since $y$ and $\epsilon$ are arbitrary, this proves that $D$ is dense in $Y$.
\end{proof}

\begin{lemma}
\label{lemma:corr}
Let $X$ be a Banach space, let $Y$ be a closed subspace and denote with $\func{\pi}{X}{X/Y}$ the quotient map. For each $r>0$, the correspondence
\[ \func{f_r}{\textup{int}(rB_{X/Y})}{2^X},\quad f_r(x)=\pi^{-1}(x)\cap rB_X \]
assumes closed, convex and nonempty values and it is lower hemicontinuous.
\end{lemma}

\begin{theo}
Let $X$ be a separable Banach space and let $\mathcal{F}$ a family of subsets of $X$ with the compactivorous property. Then there is $r>0$ with the property that, for every compact set $K\subset rB_X$, there is $x\in X$ such that $x+K\in\mathcal{F}$.
\end{theo}
\begin{proof}
By invariance under isomorphisms, we can assume that $X=\ell_1/Y$ for some closed subspace $Y\subset\ell_1$. Denote with $\func{\pi}{\ell_1}{X}$ the quotient map. By way of contradiction, suppose that for each positive integer $n$ there is $K_n\subset \textup{int}(n^{-1}B_X)$ such that no translation of $K_n$ belongs to $\mathcal{F}$. Consider for each $n$ the correspondence
\[ \func{f_n}{\textup{int}(n^{-1}B_{X})}{2^{\ell_1}},\quad f_n(x)=\pi^{-1}(x)\cap n^{-1}B_{\ell_1}, \]
which is lower hemicontinuous by Lemma \ref{lemma: corr}. Let $s_n$ be a continuous selection of $f_n$ and observe that $s_n(K_n)$ is a compact set in $n^{-1}B_{\ell_1}$. Define
\[ K=\{0\}\cup\bigcup_{n=1}^{\infty}s_n(K_n). \]
$K$ is compact and has therefore its modulus of decay $\Omega_K$. Let $L$ be the saturated compact set with the same modulus of decay. Let $V$ be an open set of $X$ such that $\pi(L)\cap V\ne\varnothing$. The set $L\cap\pi^{-1}(V)$ is relatively open in $L$, hence there are $n$ and $z\in\ell_1$ such that $z+s_n(K_n)\subset L\cap\pi^{-1}(V)$. This implies that $\pi(z)+K_n\subset \pi(L)\cap V$, thus since $\mathcal{F}$ is hereditary there cannot be $x\in X$ such that $x+\pi(L)\cap V\in\mathcal{F}$, otherwise $x+K_n\in\mathcal{F}$, against the assumption on $K_n$.
\end{proof}